\documentclass[12pt,a4paper]{article}
\usepackage{amsfonts}
\usepackage{amsmath}
\usepackage{amsthm}
\usepackage[arrow,curve,matrix]{xy}
\let\mod=\undefined
\DeclareMathOperator{\coker}{coker}
\DeclareMathOperator{\End}{End}
\DeclareMathOperator{\GL}{GL}
\DeclareMathOperator{\Grass}{Grass}
\DeclareMathOperator{\Hom}{Hom}
\DeclareMathOperator{\im}{Im}
\DeclareMathOperator{\mod}{mod}
\DeclareMathOperator{\rad}{rad}
\DeclareMathOperator{\rep}{rep}
\newcommand{\BA}{{\mathbb A}}
\newcommand{\BD}{{\mathbb D}}
\newcommand{\BM}{{\mathbb M}}
\newcommand{\CB}{{\mathcal B}}
\newcommand{\CC}{{\mathcal C}}
\newcommand{\CO}{{\mathcal O}}
\newcommand{\CU}{{\mathcal U}}
\newcommand{\CV}{{\mathcal V}}
\newcommand{\mm}{\mathfrak{m}}
\newcommand{\dd}{{\mathbf d}}
\newcommand{\ov}{\overline}
\newcommand{\bsmatrix}[1]{\left[\begin{smallmatrix} #1%
 \end{smallmatrix}\right]}
\newcommand{\psmatrix}[1]{\left(\begin{smallmatrix} #1%
 \end{smallmatrix}\right)}
\newtheorem{thm}{Theorem}[section]
\newtheorem{cor}[thm]{Corollary}
\newtheorem{lem}[thm]{Lemma}
\newtheorem{prop}[thm]{Proposition}
\newtheorem{step}{Step}
\numberwithin{equation}{section}
\begin{document}
\title{Regularity in codimension one of orbit closures
in module varieties
 }
\author{G. Zwara}
%
\maketitle

\begin{abstract}
Let $\BM_d(k)$ denote the space of $d\times d$-matrices with
coefficients in an algebraically closed field $k$.
Let $X$ be an orbit closure in the product $[\BM_d(k)]^t$ equipped
with the action of the general linear group $\GL_d(k)$ by
simultaneous conjugation.
We show that $X$ is regular at any point $y$ such that
the orbit of $y$ has codimension one in $X$.
The proof uses mainly the representation theory of associative
algebras.
\end{abstract}

\section{Introduction and the main results}

Throughout the paper, $k$ denotes an algebraically closed field
and by an algebra we mean an associative $k$-algebra with
an identity.
Let $d$ and $t$ be positive integers.
The points of $[\BM_d(k)]^t$ correspond to the algebra
homomorphisms from the free algebra
$k\langle X_1,\ldots,X_t\rangle$ to $\BM_d(k)$, or equivalently,
to the left $k\langle X_1,\ldots,X_t\rangle$-modules with
underlying vector space $k^d$.
Furthermore, the isomorphism classes of $d$-dimensional left
$k\langle X_1,\ldots,X_t\rangle$-modules correspond to
the orbits in $[\BM_d(k)]^t$ under the action of the general linear
group $\GL_d(k)$ via
$$
g\star(m_1,\ldots,m_t)=(gm_1g^{-1},\ldots,gm_tg^{-1}).
$$
Now let $A$ be a finitely generated algebra and
$a_1,\ldots,a_t$ be some generators, for a positive integer $t$.
Then we get an isomorphism
$A\simeq k\langle X_1,\ldots,X_t\rangle/I$, where $I$ is
a two-sided ideal.
Consequently, the set $\mod_A^d(k)$ of left $A$-modules with
underlying vector space $k^d$ can be identified with the
$\GL_d(k)$-invariant closed subvariety of $[\BM_d(k)]^t$
consisting of $t$-tuples $(m_1,\ldots,m_t)$ such that
$\rho(m_1,\ldots,m_t)$ is the zero matrix for any (noncommutative)
polynomial $\rho$ in $I$.
The affine variety $\mod_A^d(k)$ is called a module variety and
depends on the choice of generators of $A$ only up to
a $\GL_d(k)$-equivariant isomorphism.
We shall denote by $\CO_M$ the $\GL_d(k)$-orbit of a module $M$
in $\mod_A^d(k)$, and the closure of $\CO_M$ with respect to the
Zariski topology will be denoted by $\ov{\CO}_M$.
The main result of the paper solves the open problem posed by
Bongartz in \cite[\S 6.2,p.598]{Bmin}.

\begin{thm} \label{main1}
Let $M$ and $N$ be points in $\mod_A^d(k)$ such that $N$ belongs
to $\ov{\CO}_M$ and $\dim\CO_M-\dim\CO_N=1$.
Then the variety $\ov{\CO}_M$ is regular at $N$.
\end{thm}

Let $M$ be a module in $\mod_A^d(k)$, where $A$ is
a representation finite algebra, that is, there are only finitely
many isomorphism classes of indecomposable modules in $\mod A$.
Then $\ov{\CO}_M$ contains only finitely many $\GL_d(k)$-orbits
and hence we get the following result:

\begin{cor}
Let $A$ be a representation finite algebra and $d$ be a positive
integer. Then the closures of $\GL_d(k)$-orbits in $\mod_A^d(k)$
are regular in codimension one.
\end{cor}

We also know that such orbit closures are unibranch,
by \cite{Zuni}, but we do not know if they are normal.
The orbit closures in $\mod_A^d(k)$ are normal and Cohen-Macaulay
provided $A$ is the path algebra of a Dynkin quiver of type
$\BA_n$ or $\BD_n$ (\cite{BZ2}), or $A$ is a Brauer tree algebra
(\cite{SZder}).

We shall show in Section~\ref{proof1} that Theorem~\ref{main1}
follows from the following fact.

\begin{thm} \label{main2}
Let $0\to Z\xrightarrow{\psmatrix{\tilde{f}\\ \tilde{g}}}Z\oplus Y
\xrightarrow{(\tilde{f},-\tilde{h})}Z\to 0$ be a nonsplittable
exact sequence of finite dimensional left $A$-modules with
$Z$ indecomposable.
Then $\dim_k\End_A(Z)-\dim_k\End_A(Y)>1$.
\end{thm}

We obtain from the above exact sequence two $A$-endomorphisms
$x=\tilde{g}\tilde{h}$ and $y=\tilde{g}\tilde{f}\tilde{h}$
of the module $Y$.
These endomorphisms satisfy the relations $xy=yx$ and $x^3=y^2$,
which allows to consider $\End_A(Y)$ as a bimodule over the ring
$R=k[x,y]/(x^3-y^2)$.
Section~\ref{modbimod} is devoted to the study of properties
of modules and bimodules over the ring $R$ related to the
existence of their finite free resolutions.
Results obtained there will be used in Section~\ref{proof2} to
study the bimodule $\End_A(Y)$, leading to the proof
of Theorem~\ref{main2}.
Section~\ref{remarks} provides some consequences
of Theorem~\ref{main1} and additional remarks.

For basic background on the representation theory of algebras we
refer to \cite{ARS} and \cite{Rin}.
The author gratefully acknowledges support from the Polish
Scientific Grant KBN No.\ 5 PO3A 008 21.

\section{The proof of Theorem~\ref{main1}}
\label{proof1}

Throughout the section, $A$ is a finitely generated algebra
and $\mod A$ denotes the category of finite dimensional left
$A$-modules.
Furthermore, we abbreviate $\dim_k\Hom_A(X,Y)$ to $[X,Y]$,
for any modules $X$ and $Y$ in $\mod A$.

\begin{lem} \label{seqineq}
Let $\sigma:0\to U\xrightarrow{f}W\xrightarrow{g}V\to 0$ be
an exact sequence in $\mod A$ and $X$ be a module in $\mod A$.
Then
\begin{enumerate}
\item[(1)] $[U\oplus V,X]\geq[W,X]$ and the equality holds if and
 only if any homomorphism in $\Hom_A(U,X)$ factors through $f$;
\item[(2)] $[X,U\oplus V]\geq[X,W]$ and the equality holds if and
 only if any homomorphism in $\Hom_A(X,V)$ factors through $g$.
\end{enumerate}
\end{lem}

\begin{proof}
\textit{(1)} follows from the induced exact sequence
$$
0\to\Hom_A(V,X)\xrightarrow{\Hom_A(g,X)}\Hom_A(W,X)
\xrightarrow{\Hom_A(f,X)}\Hom_A(U,X)
$$
and \textit{(2)} follows by duality.
\end{proof}

\begin{lem} \label{equivsplit}
Let $\sigma:0\to U\xrightarrow{f}W\xrightarrow{g}V\to 0$ be
an exact sequence in $\mod A$.
Then the following conditions are equivalent:
\begin{enumerate}
\item[(1)] the sequence $\sigma$ splits;
\item[(2)] $W$ is isomorphic to $U\oplus V$;
\item[(3)] $[U\oplus V,U]=[W,U]$;
\item[(4)] $[V,U\oplus V]=[V,W]$.
\end{enumerate}
\end{lem}

\begin{proof}
Clearly the condition \textit{(1)} implies \textit{(2)}, and the
condition \textit{(2)} implies \textit{(3)} and \textit{(4)}.
Applying Lemma~\ref{seqineq} we get that \textit{(3)} implies that
the endomorphism $1_U$ factors through $f$, which means that $f$
is a section and \textit{(1)} holds.
Similarly, it follows from \textit{(4)} that $g$ is a retraction
and \textit{(1)} holds.
\end{proof}

Throughout the section, $M$ and $N$ are two modules in
$\mod_A^d(k)$ such that $N\in\ov{\CO}_M$ and
$\dim\CO_M-\dim\CO_N=1$.
Applying \cite[Theorem 1]{Zgiv} we get modules $Z$, $T$ and
the exact sequences in $\mod A$
\begin{gather}
\label{fundseq}
0\to Z\xrightarrow{f}Z\oplus M\xrightarrow{g}N\to 0,\\
\label{fundseqdual}
0\to N\xrightarrow{f'}T\oplus M\xrightarrow{g'}T\to 0.
\end{gather}

\begin{lem} \label{dimhomend}
$[M,M]=[M,N]=[N,M]=[N,N]-1$.
\end{lem}

\begin{proof}
Since the isotropy group of the point $M$ can be identified with
the automorphism group of the $A$-module $M$ and the latter is
open in the vector space $\End_A(M)$, then
$\dim\CO_M=\dim\GL_d(k)-[M,M]$.
Similarly, $\dim\CO_N=\dim\GL_d(k)-[N,N]$, which gives
$[M,M]=[N,N]-1$.
Applying Lemmas \ref{seqineq} and \ref{equivsplit} to the
sequences \eqref{fundseq} and \eqref{fundseqdual} we get the
inequalities
$$
[M,M]\leq[N,M]<[N,N]\qquad\text{and}\qquad
[M,M]\leq[M,N]<[N,N].
$$
Now the claim follows easily.
\end{proof}

Let $\rad(-,-)$ denote the two-sided ideal of the functor
$$
\Hom_A(-,-):\mod A\times\mod A\to\mod k
$$
generated by the nonisomorphisms between indecomposable modules.
From now on, we assume that $f$ belongs to $\rad(Z,Z\oplus M)$.
In fact, if this is not the case, then $f$ is of the form
$\psmatrix{f'&0\\ 0&f''}:Z'\oplus Z''\to Z'\oplus(Z''\oplus M)$,
where $f'$ is an isomorphism and $f''$ belongs to
$\rad(Z'',Z''\oplus M)$.
Consequently, the exact sequence \eqref{fundseq} has the form
$$
0\to Z'\oplus Z''\xrightarrow{\psmatrix{f'&0\\ 0&f''}}Z'\oplus
(Z''\oplus M)\xrightarrow{(0,g'')}N\to 0
$$
and we can replace it by the exact sequence
$$
0\to Z''\xrightarrow{f''}Z''\oplus M\xrightarrow{g''}N\to 0.
$$

\begin{lem} \label{getdiagram}
There is an open neighbourhood $\CU$ of $f$ in
$\Hom_A(Z,Z\oplus M)$ such that for any $f'$ in $\CU$ either
$f'$ is a section, or $f'=jfi$ for some $A$-endomorphisms $i$
and $j$ of $Z$ and $Z\oplus M$, respectively.
\end{lem}

\begin{proof}
We first recall a construction described in \cite{Zuni}
for the module $X=Z\oplus M$.
Let $c=[X,M]$.
The natural action of $\GL_d(k)$ on the space $\Hom_k(X,k^d)$
induces canonically an action of $\GL_d(k)$ on the Grassmann
variety $\Grass(\Hom_k(X,k^d),c)$ of $c$-dimensional subspaces of
the vector space $\Hom_k(X,k^d)$.
We consider the $\GL_d(k)$-variety
$$
\CC=\mod_A^d(k)\times\Grass(\Hom_k(X,k^d),c),
$$
and its one special $\GL_d(k)$-orbit
$$
\CO_{M_X}=\{(M',\Hom_A(X,M'));\;M'\in\CO_M\}.
$$
Let $\pi:\ov{\CO}_{M_X}\to\ov{\CO}_M$ denote the restriction
of the projection of $\CC$ on $\mod_A^d(k)$.

Now we want to construct a special regular morphism from
an open subset of $\Hom_A(Z,X)$ to $\ov{\CO}_{M_X}$ in a similar
way as in the proof of \cite[Proposition 3.4]{Rie}.
Let $e=\dim_kZ$.
By choosing bases, we may assume that $Z$ belongs to
$\mod_A^e(k)$ and $X$ belongs to $\mod_A^{e+d}(k)$.
Then the elements of $\Hom_A(Z,X)$ can be considered as
$(e+d)\times e$-matrices.
We choose an $(e+d)\times d$-matrix $b$ such that the matrix
$(f,b)$ is invertible.
Observe that $\dim_k\coker(\Hom_A(X,f'))=c$ for any injective
homomorphism $f':Z\to Z\oplus M$.
Let $w_1,\ldots,w_c$ be elements
of $\End_A(X)\subseteq\BM_{e+d}(k)$ whose residue classes form
a basis of $\coker(\Hom_A(X,f))$.
It is easy to see that there is an open neighbourhood $\CV$ of $f$
in $\Hom_A(Z,X)$ such that the matrix $[f',b]$ is invertible
(in particular $f'$ is injective) and the residue classes
of $w_1,\ldots,w_c$ form a basis of $\coker(\Hom_A(X,f'))$,
for any homomorphism $f'\in\CV$.
Let $f'\in\CV$, $g=[f',b]$ and write $g^{-1}=\bsmatrix{g'\\ g''}$,
where $g'$ consists of the first $e$-rows of $g^{-1}$.
Then $g^{-1}\star X=\bsmatrix{Z&W\\ 0&N'}$, that is, $N'$ is
a module in $\mod_A^d(k)$ and
$$
0\to Z\xrightarrow{f'}X\xrightarrow{g''}N'\to 0
$$
is an exact sequence in $\mod A$.
We conclude from the induced exact sequence
$$
0\to\Hom_A(X,Z)\xrightarrow{\Hom_A(X,f')}\Hom_A(X,X)\xrightarrow
{\Hom_A(X,g'')}\Hom_A(X,N')
$$
that $g''(w_1),\ldots,g''(w_c)$ form a basis of the image
$\im(\Hom_A(X,g''))$.
Hence we get a regular morphism $\Theta:\CV\to\CC$ sending $f'$
to $(N',\im(\Hom_A(X,g''))$.
If $f'\in\CV$ is a section, then $N'\in\CO_M$ and
$\im(\Hom_A(X,g''))=\Hom_A(X,N')$, and consequently, $\Theta(f')$
belongs to the orbit $\CO_{M_X}$.
Since the sections in $\CV$ form an open subset of the irreducible
set $\CV$, then the image of $\Theta$ is contained in
$\ov{\CO}_{M_X}$.
On the other hand, if $f'\in\CV$ is not a section, then
$N'$ is not isomorphic to $M$, which implies that $\Theta(f')$
belongs to the boundary
$\partial\CO_{M_X}=\ov{\CO}_{M_X}\setminus\CO_{M_X}$
of $\CO_{M_X}$.
Since $\dim\partial\CO_{M_X}<\dim\CO_{M_X}=\dim\CO_M=\dim\CO_N+1$,
the inverse image $\pi^{-1}(\CO_N)$ is a finite (disjoint)
union of $\GL_d(k)$-orbits and each of them is open in
$\partial\CO_{M_X}$.
Let $\CO_1$ denote the one containing $\Theta(f)$.
Then $\CO_{M_X}\cup\CO_1$ is an open subset of $\ov{\CO}_{M_X}$,
and consequently, $\CU=\Theta^{-1}(\CO_{M_X}\cup\CO_1)$ is an open
subset of $\CV$.

Assume that $\Theta(f')=(N',\im(\Hom_A(X,g'')))$ belongs to
$\CO_1$.
Then $N'=h\star N$ and $\im(\Hom_A(X,g''))=h\star\im(\Hom_A(X,g))$
for some element $h$ in $\GL_d(k)$.
Hence $h:N\to N'$ is an $A$-isomorphism and
$\im(\Hom_A(X,g''))=\im(\Hom_A(X,hg))$.
In particular, $hg=g''j$ for some $j\in\End_A(X)$.
Thus we obtain a commutative diagram with exact rows
$$
\xymatrix{
0\ar[r]&Z\ar[r]^-{f}\ar[d]_{i'}&X\ar[r]^-{g}\ar[d]^j
 &N\ar[r]\ar[d]^h&0\\
0\ar[r]&Z\ar[r]^-{f'}&X\ar[r]^-{g''}&N'\ar[r]&0.
}
$$
Since $h$ is an isomorphism, the sequence
$$
0\to Z\xrightarrow{\psmatrix{f\\ i'}}X\oplus Z
\xrightarrow{(j,-f')}X\to 0
$$
is exact.
Then the homomorphism $\psmatrix{f\\ i'}$ is a section,
by Lemma~\ref{equivsplit}.
The same is true for the endomorphism $i'$, as $f$ belongs to
$\rad(Z,X)$.
Hence $i'$ is an isomorphism and $f'=jfi$, where $i=(i')^{-1}$.
\end{proof}

Let $\rad_A(X)$ denote the Jacobson radical of a module
$X$ in $\mod A$ and let $\rad(E)$ denote the Jacobson radical of
an algebra $E$.
In particular, $\rad(\End_A(X))=\rad(X,X)$ for any module $X$
in $\mod A$.
Moreover, if $X$ is indecomposable, then the algebra $\End_A(X)$
is local with the maximal ideal $\rad(X,X)$ consisting of the
nilpotent endomorphisms of $X$ and there is a decomposition
$\End_A(X)=k\cdot 1_X\oplus\rad(X,X)$.

\begin{lem} \label{Zind}
The module $Z$ is indecomposable, $[N,Z]=[M,Z]+1$ and any
radical endomorphism of $Z$ factors through $f$.
\end{lem}

\begin{proof}
Suppose that $Z=Z_1\oplus Z_2$ for two nonzero modules $Z_1$ and
$Z_2$.
Since $f$ belongs to $\rad(Z_1\oplus Z_2,Z_1\oplus Z_2\oplus M)$,
then the map $f+t\cdot 1_{Z_1}:Z\to Z\oplus M$ is not a section,
for any $t\in k$.
Applying Lemma~\ref{getdiagram}, we get $f+t\cdot 1_{Z_1}=jfi$
for some $t\neq 0$ and endomorphisms $i$ and $j$.
Since $f$ belongs to $\rad(Z,Z\oplus M)$, the same holds
for $f+t\cdot 1_{Z_1}$ and $t\cdot 1_{Z_1}$, a contradiction.
Therefore the module $Z$ is indecomposable.

Let $E=\End_A(Z)$.
We have the induced exact sequence in $\mod E$
$$
0\to\Hom_A(N,Z)\xrightarrow{\Hom_A(g,Z)}\Hom_A(Z\oplus M,Z)
\xrightarrow{\Hom_A(f,Z)}\Hom_A(Z,Z).
$$
Then the image of $\alpha=\Hom_A(f,Z)$ is contained in
$\rad(E)=\rad(Z,Z)$ as $f$ belongs to $\rad(Z,Z\oplus M)$.
It remains to show the reverse inclusion, which means that
the restriction
$$
\alpha':\Hom_A(Z\oplus M,Z)\to\rad(E)
$$
of $\alpha$ is surjective.
Since $\im(\alpha')$ is an $E$-submodule and
$\rad_E(\rad(E))=\rad^2(E)$, it suffices to show that
the composition
$$
\beta:\Hom_A(Z\oplus M,Z)\to\rad(E)/\rad^2(E)
$$
of $\alpha'$ followed by a quotient is surjective.

Let $h\in\rad(E)$.
Observe that $f+t\cdot\psmatrix{h\\ 0}$ belongs to
$\rad(Z,Z\oplus M)$ for any $t\in k$.
Applying Lemma~\ref{getdiagram} we get $f+t\cdot\psmatrix{h\\ 0}=
jfi$ for some $t\neq 0$ and endomorphisms $i$ and $j$.
Then we have the equality $(1_Z,0)f+t\cdot h=j'fi$ in $E$, where
$j'=(1_Z,0)j:Z\oplus M\to Z$.
We decompose $i=c\cdot 1_Z+i'$, where $c\in k$ and $i'\in\rad(E)$.
Since $j'f$ belongs to $\rad(E)$, then
$j'fi-cj'f$ belongs to $\rad^2(E)$.
Altogether, we conclude that
$$
h+\rad^2(E)=t^{-1}\cdot(cj'-(1_Z,0))f+\rad^2(E)
=\beta(ct^{-1}j'-(t^{-1}\cdot 1_Z,0)),
$$
which finishes the proof.
\end{proof}

We decompose $f=\psmatrix{f_1\\ f_2}$ and $g=(g_1,g_2)$.
Then the square
\begin{equation} \label{funddiag}
\xymatrix{
Z\ar[r]^{f_2}\ar[d]_{f_1}&M\ar[d]^{-g_2}\\ Z\ar[r]^{g_1}&N
}
\end{equation}
is exact, that is, it is a pushout and a pull-back.
Furthermore $f_1$ is nilpotent.

\begin{lem} \label{powerofM}
Let $j$ be a positive integer such that $(f_1)^j=0$.
Then any radical endomorphism of $Z$ factors through
$(f_1,b):Z\oplus M^j\to Z$ for some $A$-homomorphism $b$.
\end{lem}

\begin{proof}
Let $e'$ be an element of $\End_A(Z)$.
Applying Lemma~\ref{Zind} we obtain a decomposition
$$
e'=\lambda\cdot 1_Z+a'f_1+b'f_2
$$
for some scalar $\lambda\in k$ and $A$-homomorphisms $a':Z\to Z$
and $b':M\to Z$.
Let $e$ be a radical endomorphism of $Z$.
Using the above $j$ times, we get
$$
e=\sum_{i=0}^{j-1}\lambda_i\cdot(f_1)^i+b_if_2(f_1)^i.
$$
Since the endomorphism $e$ is radical then $\lambda_0=0$.
Hence we get the claim for $b=(b_0,\ldots,b_{j-1})$.
\end{proof}

\begin{prop} \label{longprop}
Assume that Theorem~\ref{main2} holds. Then $[Z,M]=[Z,N]$.
\end{prop}

\noindent\textit{Proof of Proposition~\ref{longprop}.}
Suppose that $[Z,M]\neq[Z,N]$.
We divide the proof into several steps.

\begin{step} \label{step1}
$g_1$ factors through $f=\psmatrix{f_1\\ f_2}$.
\end{step}

\begin{proof}
Applying Lemma~\ref{seqineq} for $X=Z$ and the sequence
\eqref{fundseq} we get a homomorphism $u$ in $\Hom_A(Z,N)$
which does not factor through $g$.
Furthermore, we may assume that $uf_1$ factors through $g$
as $f_1$ is nilpotent.
Hence
\begin{equation} \label{x1}
uf_1=g_1a_1+g_2a_2
\end{equation}
for some $A$-homomorphisms $a_1:Z\to Z$ and $a_2:Z\to M$.
The homomorphism $a_2$ factors through $f$, by
Lemma~\ref{dimhomend} and Lemma~\ref{seqineq} applied for $X=M$
and the sequence \eqref{fundseq}.
We decompose the endomorphism $a_1=\lambda\cdot 1_Z+a'_1$,
where $\lambda\in k$ and $a'_1$ belongs to $\rad(Z,Z)$.
By Lemma~\ref{Zind}, $a'_1$ also factors through $f$, and
consequently,
$$
a_1=\lambda\cdot 1_Z+b_{1,1}f_1+b_{1,2}f_2,\qquad
a_2=b_{2,1}f_1+b_{2,2}f_2
$$
for some $A$-homomorphisms $b_{1,1}$, $b_{1,2}$, $b_{2,1}$ and
$b_{2,2}$.
Combining these equalities with \eqref{x1} we get
\begin{equation} \label{x2}
\lambda\cdot g_1=(u-g_1b_{1,1}-g_2b_{2,1})f_1+(-g_1b_{1,2}-g_2
b_{2,2})f_2.
\end{equation}

Suppose that $\lambda=0$.
Then it follows from the exactness of \eqref{funddiag} that
$$
u-g_1b_{1,1}-g_2b_{2,1}=cg_1
$$
for some $A$-endomorphism $c:N\to N$.
We know that $[N,N]-[N,M]=1$, by Lemma~\ref{dimhomend}.
We conclude from the induced exact sequence
$$
0\to\Hom_A(N,Z)\xrightarrow{\Hom_A(N,f)}\Hom_A(N,Z\oplus M)
\xrightarrow{\Hom_A(N,g)}\Hom_A(N,N)
$$
that $\Hom_A(N,N)=k\cdot 1_N\oplus\im(\Hom_A(N,g))$.
Hence $c=\mu\cdot 1_N+g_1d_1+g_2d_2$ for some
$\mu\in k$ and $A$-homomorphisms $d_1$ and $d_2$.
Consequently, the homomorphism
$$
u=g_1(b_{1,1}+\mu\cdot 1_Z+d_1g_1)+g_2(b_{2,1}+d_2g_1)
$$
factors through $g=(g_1,g_2)$, a contradiction.
Thus $\lambda\neq 0$ and the equality \eqref{x2} shows that $g_1$
factors through $f=\psmatrix{f_1\\ f_2}$.
\end{proof}

Hence the exact square \eqref{funddiag} divides into two exact
squares
\begin{equation} \label{fund2diag}
\xymatrixcolsep{3pc}
\xymatrix{
Z\ar[r]^-{u}\ar[d]_{f_1}&W\ar[r]^-{x}\ar[d]^{\psmatrix{w_1\\ w_2}}
 &M\ar[d]^{-g_2}\\
Z\ar[r]^-{\psmatrix{f_1\\ f_2}}&Z\oplus M\ar[r]^-{(y_1,y_2)}&N.
}
\end{equation}

\begin{step} \label{step2}
The homomorphism $w_2:W\to M$ is a retraction and the inequality
$[Z,Z\oplus M]-[Z,W]\leq 1$ holds.
\end{step}

\begin{proof}
Applying Lemma~\ref{seqineq} for $X=M$ and the exact squares
\eqref{fund2diag}, we get that the integers
$$
[M,Z\oplus M]-[M,W]\qquad\text{and}\qquad
[M,W\oplus N]-[M,M^2\oplus Z]
$$
are nonnegative.
Moreover, their sum equals $[M,N]-[M,M]=0$,
by Lemma~\ref{dimhomend}.
Hence these numbers are zero and any map in $\Hom_A(M,Z\oplus M)$
factors through $\psmatrix{f_1&w_1\\ f_2&w_2}$, by
Lemma~\ref{seqineq} applied for $X=M$ and the left square in
\eqref{fund2diag}.
Consequently, any map in $\Hom_A(M,Z)$ factors through $(f_1,w_1)$
while any endomorphism in $\End_A(M)$ factors through $(f_2,w_2)$.
In particular, $(f_2,w_2)$ is a retraction and the same holds for
$w_2$, as $f_2$ belongs to $\rad(Z,M)$.
Furthermore, $\Hom_A(Z,M)$ is contained in the image
of the map $\alpha=\Hom_A(Z,\psmatrix{f_1&w_1\\ f_2&w_2})$
in the induced exact sequence
$$
0\to\Hom_A(Z,Z)\to\Hom_A(Z,Z\oplus W)\xrightarrow{\alpha}
\Hom_A(Z,Z\oplus M).
$$
Hence the inequality $[Z,Z\oplus M]-[Z,W]\leq 1$ will be
a consequence of the fact that $\rad(Z,Z)$ is contained in the
image of the map
$$
\Hom_A(Z,(f_1,w_1)):\Hom_A(Z,Z\oplus W)\to\Hom_A(Z,Z).
$$
The latter follows from Lemma~\ref{powerofM} and the fact that any
homomorphism in $\Hom_A(M,Z\oplus M)$ factors through $(f_1,w_1)$.
\end{proof}

Consequently, we can decompose $W=Y\oplus M$ for some $A$-module
$Y$, such that $w_2=(0,1_M):Y\oplus M\to M$.

\begin{step} \label{step3}
There is a nonsplittable exact sequence in $\mod A$ of the form
\begin{equation} \label{firstseq}
0\to Z\xrightarrow{\psmatrix{\tilde{f}\\ \tilde{g}}}Z\oplus Y
\xrightarrow{(\tilde{f},-\tilde{h})}Z\to 0.
\end{equation}
\end{step}

\begin{proof}
We decompose $u=\psmatrix{u_1\\ u_2}:Z\to Y\oplus M$ and
$w_1=(v_1,v_2):Y\oplus M\to Z$.
We conclude from \eqref{fund2diag} the exactness of the upper row
in the diagram
\begin{equation} \label{diag3}
\xymatrixcolsep{2.5pc}
\xymatrixrowsep{3pc}
\xymatrix{
0\ar[r]&Z\ar[rr]^-{\psmatrix{f_1\\ u_1\\ u_2}}\ar@{=}[d]
 &&Z\oplus Y\oplus M\ar[rr]^-{\psmatrix{f_1&-v_1&-v_2\\ f_2&0&-1}}
 \ar[d]_{\psmatrix{1&0&0\\ 0&1&0\\ -f_2&0&1}}&&Z\oplus M\ar[r]
 \ar[d]^{\psmatrix{1&-v_2\\ 0&-1}}&0\\
0\ar[r]&Z\ar[rr]^-{\psmatrix{f_1\\ u_1\\ 0}}&&Z\oplus Y\oplus M
 \ar[rr]^-{\psmatrix{f_1-v_2f_2&-v_1&0\\ 0&0&1}}&&Z\oplus M\ar[r]
 &0.
}
\end{equation}
In particular, $f_2f_1-u_2=0$, which implies that the diagram
\eqref{diag3} is commutative.
Since the maps corresponding to vertical arrows are isomorphisms,
the bottom row is exact as well.
It follows from the construction of the squares \eqref{fund2diag}
that $f_2=xu$.
We decompose $x=(x_1,x_2):Y\oplus M\to M$.
Then
$$
f_2=x_1u_1+x_2u_2=x_1u_1+x_2f_2f_1,
$$
and consequently,
$$
f_1-v_2f_2=(1_Z-v_2x_2f_2)f_1-v_2x_1u_1.
$$
Since $f_2$ belongs to $\rad(Z,M)$, the endomorphism
$a=1_Z-v_2x_2f_2$ is an isomorphism.
Then
$$
f_1-v_2f_2=af_1+abu_1,
$$
where $b=-a^{-1}v_2x_1$.
The exactness of the bottom row in the diagram \eqref{diag3}
implies the exactness of the upper row in the commutative
diagram
$$
\xymatrixcolsep{2.5pc}
\xymatrixrowsep{3pc}
\xymatrix{
0\ar[r]&Z\ar[rr]^-{\psmatrix{f_1\\ u_1}}\ar@{=}[d]&&Z\oplus Y
 \ar[rr]^-{(af_1+abu_1,-v_1)}\ar[d]_{\psmatrix{1&b\\ 0&1}}&&
 Z\ar[r]\ar[d]^{(a^{-1})}&0\\
0\ar[r]&Z\ar[rr]^-{\psmatrix{f_1+bu_1\\ u_1}}&&Z\oplus Y
 \ar[rr]^-{(f_1+bu_1,-\tilde{h})}&&Z\ar[r]&0,
}
$$
where $\tilde{h}=a^{-1}v_1+f_1b+bu_1b$.
Since the maps corresponding to the vertical arrows are
isomorphisms, the bottom row is also exact.
Setting $\tilde{f}=f_1+bu_1$ and $\tilde{g}=u_1$ we get the exact
sequence \eqref{firstseq}.
Suppose that the sequence \eqref{firstseq} splits.
Since $\tilde{f}=a^{-1}f_1-a^{-1}v_2f_2$ belongs to $\rad(Z,Z)$,
then $\tilde{g}:Z\to Y$ is a section and $\tilde{h}$ is
a retraction.
Hence both of them are isomorphisms as $\dim_kY=\dim_kZ$.
Consequently, $\tilde{f}\tilde{f}=\tilde{h}\tilde{g}$ is
an isomorphism, a contradiction.
Therefore the exact sequence \eqref{firstseq} does not split.
\end{proof}

The right square in \eqref{fund2diag} leads to the exact sequence
$$
0\to Y\oplus M\xrightarrow{\psmatrix{x_1&x_2\\ v_1&v_2\\ 0&1}}
M\oplus Z\oplus M\xrightarrow{(g_2,y_1,y_2)}N\to 0,
$$
which implies the exactness of the sequence
\begin{equation} \label{secondseq}
0\to Y\xrightarrow{\psmatrix{x_1\\ v_1}}M\oplus Z
\xrightarrow{(g_2,y_1)}N\to 0.
\end{equation}

\begin{step} \label{step4}
$[Z,Z]-[Y,Y]=1$.
\end{step}

\begin{proof}
We claim that $[Z,Y]=[Y,Y]$.
Assume first that the exact sequence \eqref{secondseq} splits.
Since the sequence \eqref{firstseq} does not split, then
$Y$ is not isomorphic to $Z$.
Hence $v_1$ belongs to $\rad(Y,Z)$, as $Z$ is indecomposable
and $\dim_kZ=\dim_kY$.
This implies that $x_1:Y\to M$ is a section.
In particular, $M$ is isomorphic to $Y\oplus Y'$ for some
$A$-module $Y'$.
Applying Lemma~\ref{seqineq} to \eqref{firstseq} we get
$[Y,Y]\leq[Z,Y]$ and $[Y,Y']\leq[Z,Y']$, and applying it to
\eqref{secondseq} we get $[M\oplus Z,M]\leq[Y\oplus N,M]$.
Consequently,
$$
0\leq[Z,Y]-[Y,Y]\leq[Z,M]-[Y,M]\leq[N,M]-[M,M]=0,
$$
by Lemma~\ref{dimhomend}.

Assume now that the exact sequence \eqref{secondseq} does not split.
Then
\begin{align*}
[M,Y\oplus N]&\geq[M,M\oplus Z],&[N,Z]&\geq[N,Y],\\
[Y\oplus N,Y]&>[M\oplus Z,Y],&[Z,Y]&\geq[Y,Y],
\end{align*}
by Lemmas \ref{seqineq} and \ref{equivsplit} applied to
the sequences \eqref{firstseq} and \eqref{secondseq}.
From Lemmas \ref{dimhomend} and \ref{Zind} we get
 $[M,M]=[M,N]$, $[N,Z]-[M,Z]=1$ and hence
\begin{align*}
0&\leq[Z,Y]-[Y,Y]\leq[N,Y]-[M,Y]-1\\
&\leq([N,Z]-[M,Z]-1)+([M,Z]-[M,Y])\leq[M,N]-[M,M]=0,
\end{align*}
which proves the claim.

By Step~\ref{step2}, we get $[Z,Z]-[Z,Y]\leq 1$.
But $[Z,Z]>[Z,Y]$, by Step~\ref{step3} and Lemmas \ref{seqineq}
and \ref{equivsplit}.
Therefore
$
[Z,Z]-[Y,Y]=[Z,Z]-[Z,Y]=1.
$
\end{proof}

Steps \ref{step3} and \ref{step4} give a contradiction with
Theorem~\ref{main2}.
This finishes the proof of Proposition~\ref{longprop}.
\qed
\bigskip

\noindent\textit{Deduction of Theorem~\ref{main1} from
Theorem~\ref{main2}.}
Applying Lemma~\ref{dimhomend} and Proposition~\ref{longprop}
we get $[Z\oplus M,M]=[Z\oplus M,N]$.
Then the variety $\ov{\CO}_M$ is regular at the point $N$,
by \cite[Proposition 2.2]{Zuni}.
\qed

\section{Bimodules over $k[x,y]/(x^3-y^2)$}
\label{modbimod}

Let $R=k[m^2,m^3]$ denote the subalgebra of the polynomial ring
$k[m]$ in a formal variable $m$.
We say that a left $R$-module $M$ has property [P1] if the
sequence
$$
\begin{pmatrix}M\\ M\end{pmatrix}
\xrightarrow{\psmatrix{m^3&-m^2\\ m^4&-m^3}\cdot}
\begin{pmatrix}M\\ M\end{pmatrix}
\xrightarrow{\psmatrix{m^3&-m^2\\ m^4&-m^3}\cdot}
\begin{pmatrix}M\\ M\end{pmatrix}
$$
is exact.
We shall see later (Corollary~\ref{finresmod}) that this is
equivalent to the fact that $M$ has a free resolution of finite
length.
Dually, we say that a right $R$-module $M$ has property [P1'] if
the sequence
$$
\begin{pmatrix}M& M\end{pmatrix}
\xrightarrow{\cdot\psmatrix{m^3&m^4\\ -m^2&-m^3}}
\begin{pmatrix}M& M\end{pmatrix}
\xrightarrow{\cdot\psmatrix{m^3&m^4\\ -m^2&-m^3}}
\begin{pmatrix}M& M\end{pmatrix}
$$
is exact.
Observe that $\mm=(m^2,m^3)$ is a maximal ideal of $R$.

\begin{lem} \label{subfree}
Let $M$ be a submodule of a left free $R$-module.
If $M$ has property [P1] then it is free.
\end{lem}

\begin{proof}
Let $\{b_s\}_{s\in S}$ be a set of elements of $M$ whose residue
classes form a linear basis of $M/\mm M$.
We want to show that this set is a basis of the $R$-module $M$.
Since $M$ is contained in a free $R$-module $W$ then
$$
\bigcap_{i\geq 1}\mm^i M\subseteq\bigcap_{i\geq 1}\mm^i W=\{0\}.
$$
By Nakayama's lemma, the elements $b_s$, $s\in S$ generate the
$R$-module $M$.

Assume that $\sum_{s\in S}r_s b_s=0$, where all but a finite
number of elements $r_s\in R$ are zero.
We decompose $r_s=a_{s,0}+\sum_{i\geq 2}a_{s,i}m^i$, $s\in S$,
where $a_{s,i}$ are scalars in $k$.
It follows from the definition of $b_s$, $s\in S$ that
$a_{s,0}=0$ for any $s\in S$.
We have to show that $r_s=0$ for any $s\in S$, which means that
$a_{s,i}=0$ for any $s\in S$ and $i\geq 2$.
Suppose this is not the case and let $j\geq 2$ denote the minimal
integer such that there is some $s_0\in S$ with $a_{s_0,j}\neq 0$.
Then
$$
0=m^3(\sum_{s\in S}(\sum_{i\geq j}a_{s,i}m^i)b_s)=
m^j(\sum_{s\in S}(\sum_{i\geq 3}a_{s,j+i-3}m^i)b_s).
$$
Since $M$ is contained in a free $R$-module, $m^j$ is not a zero
divisor in $M$.
Consequently,
$$
\sum_{s\in S}(\sum_{i\geq 3}a_{s,j+i-3}m^i)b_s=0.
$$
Then $m^3x'-m^2x''=0$ for
$$
x'=\sum_{s\in S}a_{s,j}b_s\qquad\text{and}\qquad
x''=-\sum_{s\in S}(\sum_{i\geq 2}a_{s,j+i-1}m^i)b_s.
$$
Moreover,
$$
0=m^3(m^3x'-m^2x'')=m^2(m^4x'-m^3x'').
$$
Since $m^2$ is not a zero divisor in $M$ and $M$ has property
[P1], we get
$$
\psmatrix{m^3&-m^2\\  m^4&-m^3}\psmatrix{x'\\ x''}
 =\psmatrix{0\\ 0}\qquad\text{and}\qquad
\psmatrix{m^3&-m^2\\  m^4&-m^3}\psmatrix{y'\\ y''}
 =\psmatrix{x'\\ x''}
$$
for some $y',y''\in M$.
Therefore $x'$ belongs to $\mm M$, and consequently, $a_{s,j}=0$
for any $s\in S$.
This gives a contradiction with the choice of $j$ and hence
the module $M$ is free.
\end{proof}

Let $M$ be an $R$-$R$-bimodule, $N=\begin{pmatrix}M&M\\ M&M
\end{pmatrix}$ and consider the maps
$$
\xi:N\xrightarrow{\psmatrix{m^3&-m^2\\ m^4&-m^3}\cdot}N
\qquad\text{and}\qquad
\eta:N\xrightarrow{\cdot\psmatrix{m^3&m^4\\ -m^2&-m^3}}N,
$$
given by left and right, respectively, multiplications of $N$ by
$2\times 2$-matrices.
Observe that $\xi\xi=\eta\eta=0$ and $\xi\eta=\eta\xi$.
We say that $M$ has property [P2] if the sequence
$$
N\xrightarrow{\xi\eta}N\xrightarrow{\psmatrix{\xi\\ \eta}}
N\oplus N
$$
is exact.
In fact, we shall see (Corollary~\ref{finresbimod}) that property
[P2] is equivalent to the fact that the bimodule $M$ has a free
resolution of finite length.

\begin{lem} \label{bigM}
Let $N=\begin{pmatrix}M&M\\ M&M\end{pmatrix}$, where $M$ is an
$R$-$R$-bimodule having property [P2]. Then $M$ has properties
[P1] and [P1'], and the following sequence is exact:
\begin{equation} \label{longN}
N\oplus N\xrightarrow{\psmatrix{\xi&\eta}}N\xrightarrow{\xi\eta}
N\xrightarrow{\psmatrix{\xi\\ \eta}}N\oplus N\xrightarrow{
\psmatrix{\xi&0\\ \eta&-\xi\\ 0&\eta}}N\oplus N\oplus N.
\end{equation}
\end{lem}

\begin{proof}
Observe that $M$ has property [P1] if and only if the sequence
$$
N\xrightarrow{\xi}N\xrightarrow{\xi}N
$$
is exact.
We take $n\in N$ such that $\xi(n)=0$ and set $n_1=\eta(n)$.
Then $\xi(n_1)=\eta(n_1)=0$, which implies that $n_1=\eta\xi(n_2)$
for some $n_2\in N$.
Let $n_3=n-\xi(n_2)$.
Then $\xi(n_3)=\eta(n_3)=0$, which gives $n_3=\xi\eta(n_4)$ for
some $n_4\in N$.
Consequently, $n=\xi(n_2+\eta(n_4))$.
This shows that the sequence $N\xrightarrow{\xi}N\xrightarrow{\xi}
N$ is exact.
By a similar diagram chasing, one can get the exactness of
the sequences $N\xrightarrow{\eta}N\xrightarrow{\eta}N$ and
\eqref{longN}, which proves the claim.
\end{proof}

\begin{lem} \label{transitive}
Let $0\to M_1\to M_2\to M_3\to 0$ be an exact sequence
of $R$-$R$-bimodules.
If two of the bimodules $M_1$, $M_2$ and $M_3$ have property [P2],
then the third one has as well.
\end{lem}

\begin{proof}
Assume that $0\to M_1\to M_2\to M_3\to 0$ is an exact sequence
of $R$-$R$-bimodules.
Then we get the exact sequence $0\to N_1\to N_2\to N_3\to 0$,
where $N_i=\begin{pmatrix}M_i&M_i\\ M_i&M_i\end{pmatrix}$ for
$i=1,2,3$.
We apply Lemma~\ref{bigM} and consider the commutative diagram
with exact rows
$$
\xymatrix{
&&(N_2)^2\ar[r]\ar[d]^{\psmatrix{\xi&\eta}}&(N_3)^2\ar[r]
 \ar[d]^{\psmatrix{\xi&\eta}}&0\\
0\ar[r]&N_1\ar[r]\ar[d]^{\xi\eta}&N_2\ar[r]\ar[d]^{\xi\eta}
 &N_3\ar[r]\ar[d]^{\xi\eta}&0\\
0\ar[r]&N_1\ar[r]\ar[d]^{\psmatrix{\xi\\ \eta}}&N_2\ar[r]
 \ar[d]^{\psmatrix{\xi\\ \eta}}&N_3\ar[r]\ar[d]^{\psmatrix{\xi\\
 \eta}}&0\\
0\ar[r]&(N_1)^2\ar[r]\ar[d]^{\psmatrix{\xi&0\\ \eta&-\xi\\ 0&\eta}}
 &(N_2)^2\ar[r]\ar[d]^{\psmatrix{\xi&0\\ \eta&-\xi\\ 0&\eta}}
 &(N_3)^2\ar[r]&0\\
0\ar[r]&(N_1)^3\ar[r]&(N_2)^3.
}
$$
If two of the bimodules $M_1$, $M_2$ and $M_3$ have property [P2]
then the corresponding two columns are exact.
Hence we get the exactness in the middle of the third column,
which means that the third bimodule has also property [P2].
\end{proof}

\begin{lem} \label{freehasP2}
Any free $R$-$R$-bimodule has property [P2].
\end{lem}

\begin{proof}
Let $M$ be a free $R$-$R$-bimodule and choose a basis
$\{b_s\}_{s\in S}$.
Assume that $x_{1,1}$, $x_{1,2}$, $x_{2,1}$ and $x_{2,2}$ are
elements in $M$ such that
\begin{equation} \label{checkforfree}
\psmatrix{m^3&-m^2\\  m^4&-m^3}\psmatrix{x_{1,1}&x_{1,2}\\ x_{2,1}
&x_{2,2}}=\psmatrix{0&0\\ 0&0}\quad\text{and}\quad
\psmatrix{x_{1,1}&x_{1,2}\\ x_{2,1}&x_{2,2}}\psmatrix{m^3&m^4\\
-m^2&-m^3}=\psmatrix{0&0\\ 0&0}.
\end{equation}
We decompose
$$
x_{p,q}=\sum_{s\in S}\sum_{\substack{i\geq 0\\ i\neq 1}}
\sum_{\substack{j\geq 0\\ j\neq 1}}a^{p,q}_{s,i,j}m^ib_sm^j,
\qquad p,q=1,2,
$$
where all but a finite number of scalars $a^{p,q}_{s,i,j}$ in $k$
are zero.
We conclude from \eqref{checkforfree} that
$$
a^{1,1}_{s,i,j}=a^{1,2}_{s,i,j+1}=a^{2,1}_{s,i+1,j}=
a^{2,2}_{s,i+1,j+1}\qquad\text{for }s\in S,\;i,j\geq 2,
$$
and the remaining scalars $a^{p,q}_{s,i,j}$ are zero.
If we set
\begin{align*}
y_{1,1}&=\sum_{s\in S}a^{1,1}_{s,3,3}b_s,
&y_{1,2}&=-\sum_{s\in S}\sum_{\substack{j\geq 0\\ j\neq 1}}
 a^{1,1}_{s,3,j+2}b_sm^j,\\
y_{2,1}&=-\sum_{s\in S}\sum_{\substack{i\geq 0\\ i\neq 1}}
 a^{1,1}_{s,i+2,3}m^ib_s,
&y_{2,2}&=\sum_{s\in S}\sum_{\substack{i\geq 0\\ i\neq 1}}
 \sum_{\substack{j\geq 0\\ j\neq 1}}a^{1,1}_{s,i+2,j+2}m^ib_sm^j,
\end{align*}
then
$$
\psmatrix{x_{1,1}&x_{1,2}\\ x_{2,1}&x_{2,2}}=
\psmatrix{m^3&-m^2\\ m^4&-m^3}\psmatrix{y_{1,1}&y_{1,2}\\
y_{2,1}&y_{2,2}}\psmatrix{m^3&m^4\\ -m^2&-m^3}.
$$
Hence $M$ has property [P2].
\end{proof}

\begin{lem} \label{PLRtoPL}
Assume that $M$ is an $R$-$R$-bimodule having property [P2] which
is torsion free as a right $R$-module.
Then the left $R$-module $M/M\mm$ has property [P1] and the
following sequence is exact:
\begin{equation} \label{PRspecial}
M/Mm^2\xrightarrow{\cdot m^3}M/Mm^2\xrightarrow{\cdot m^3}M/Mm^2.
\end{equation}
\end{lem}

\begin{proof}
Since $m^2$ is not a zero divisor of the right $R$-module $M$,
then the sequence
$$
0\to M\xrightarrow{\cdot m^2}M\to M/Mm^2\to 0
$$
is exact.
In fact, this is a sequence of $R$-$R$-bimodules since the algebra
$R$ is commutative.
Then $M/Mm^2$ has property [P2], by Lemma~\ref{transitive}.
Applying Lemma~\ref{bigM} we get the exact sequence
$$
N\oplus N\xrightarrow{\psmatrix{\xi&\eta}}N\xrightarrow{\xi\eta}
N,\qquad\text{where }
N=\begin{pmatrix}M/Mm^2&M/Mm^2\\ M/Mm^2&M/Mm^2\end{pmatrix}.
$$
We have to show the exactness of the sequence
$$
\begin{pmatrix}M/M\mm\\ M/M\mm\end{pmatrix}
\xrightarrow{\psmatrix{m^3&-m^2\\  m^4&-m^3}\cdot}
\begin{pmatrix}M/M\mm\\ M/M\mm\end{pmatrix}
\xrightarrow{\psmatrix{m^3&-m^2\\  m^4&-m^3}\cdot}
\begin{pmatrix}M/M\mm\\ M/M\mm\end{pmatrix}.
$$
Let $x_1$ and $x_2$ be elements in $M$ such that
$$
\psmatrix{m^3&-m^2\\  m^4&-m^3}\begin{pmatrix}x_1+M\mm\\ x_2+M\mm
\end{pmatrix}=\begin{pmatrix}0+M\mm\\ 0+M\mm\end{pmatrix}.
$$
Then
$$
\psmatrix{m^3&-m^2\\  m^4&-m^3}\begin{pmatrix}x_1&0\\ x_2&0
 \end{pmatrix}\psmatrix{m^3&m^4\\ -m^2&-m^3}\in
\begin{pmatrix}M\mm&0\\ M\mm&0\end{pmatrix}\psmatrix{m^3&m^4\\
 -m^2&-m^3}\subseteq
\begin{pmatrix}Mm^2&Mm^2\\ Mm^2&Mm^2\end{pmatrix},
$$
hence
$$
\xi\eta
\begin{pmatrix}x_1+Mm^2&0+Mm^2\\ x_2+Mm^2&0+Mm^2\end{pmatrix}
=\begin{pmatrix}0+Mm^2&0+Mm^2\\ 0+Mm^2&0+Mm^2\end{pmatrix},
$$
and consequently,
\begin{align*}
\begin{pmatrix}x_1+Mm^2&0+Mm^2\\ x_2+Mm^2&0+Mm^2\end{pmatrix}=
&\psmatrix{m^3&-m^2\\  m^4&-m^3}\begin{pmatrix}y_1+Mm^2&y_2+Mm^2\\
 y_3+Mm^2&y_4+Mm^2\end{pmatrix}\\
&+\begin{pmatrix}y_5+Mm^2&y_6+Mm^2\\ y_7+Mm^2&y_8+Mm^2\end{pmatrix}
 \psmatrix{m^3&m^4\\ -m^2&-m^3}
\end{align*}
for some $y_1,\ldots,y_8\in M$.
This implies that
$$
\begin{pmatrix}x_1+M\mm\\ x_2+M\mm\end{pmatrix}=
\psmatrix{m^3&-m^2\\  m^4&-m^3}
\begin{pmatrix}y_1+M\mm\\ y_3+M\mm\end{pmatrix}.
$$
Therefore $M/M\mm$ has property [P1].
We know that $M/Mm^2$ has property [P1'], by Lemma~\ref{bigM}.
This gives the exact sequence
\begin{multline*}
\begin{pmatrix}M/Mm^2&M/Mm^2\end{pmatrix}\xrightarrow{\cdot
 \psmatrix{m^3&0\\ 0&-m^3}}\begin{pmatrix}M/Mm^2&M/Mm^2
\end{pmatrix}\xrightarrow{\cdot\psmatrix{m^3&0\\ 0&-m^3}}\\
\xrightarrow{\cdot\psmatrix{m^3&0\\ 0&-m^3}}
\begin{pmatrix}M/Mm^2&M/Mm^2\end{pmatrix},
\end{multline*}
from which we derive the exactness of \eqref{PRspecial}.
\end{proof}

\begin{prop} \label{freeres}
Let $M$ be an $R$-$R$-bimodule having property [P2].
If $M$ is torsion free as a right $R$-module, then there is a free
bimodule resolution
$$
0\to U\to W\to M\to 0.
$$
Furthermore, if the bimodule $M$ is finitely generated then we may
assume the same for $U$ and $W$.
\end{prop}

\begin{proof}
We take an exact sequence of $R$-$R$-bimodules
$$
0\to U\to W\to M\to 0
$$
such that the bimodule $W$ is free.
Obviously $W$ can be finitely generated provided $M$ is finitely
generated.
Since the $R$-$R$-bimodules can be equivalently considered as
$R\otimes R$-modules and the ring $R\otimes R$ is noetherian, then
the bimodule $U$ is finitely generated if $W$ is.
Furthermore, $U$, $W$ and $M$ are torsion free right $R$-modules.
Hence we get the following commutative diagram with exact columns
and upper two rows:
$$
\xymatrix{
&0\ar[d]&0\ar[d]&0\ar[d]\\
0\ar[r]&U\ar[r]\ar[d]^{\cdot m^2}&W\ar[r]\ar[d]^{\cdot m^2}&
 M\ar[r]\ar[d]^{\cdot m^2}&0\\
0\ar[r]&U\ar[r]\ar[d]&W\ar[r]\ar[d]&M\ar[r]\ar[d]&0\\
0\ar[r]&U/Um^2\ar[r]\ar[d]&W/Wm^2\ar[r]\ar[d]&M/Mm^2\ar[r]\ar[d]
 &0\\
&0&0&0.
}
$$
Consequently, the bottom row is also exact.
Applying Lemma \eqref{PRspecial} we get another commutative
diagram with exact columns and upper three rows
$$
\xymatrix{
&&W/Wm^2\ar[r]\ar[d]^{\cdot m^3}&M/Mm^2\ar[r]\ar[d]^{\cdot m^3}
 &0\\
0\ar[r]&U/Um^2\ar[r]\ar[d]^{\cdot m^3}&W/Wm^2\ar[r]
 \ar[d]^{\cdot m^3}&M/Mm^2\ar[r]\ar[d]^{\cdot m^3}&0\\
0\ar[r]&U/Um^2\ar[r]\ar[d]&W/Wm^2\ar[r]\ar[d]&M/Mm^2\ar[r]\ar[d]
 &0\\
0\ar[r]&U/U\mm\ar[r]\ar[d]&W/W\mm\ar[r]\ar[d]&M/M\mm\ar[r]\ar[d]
 &0\\
&0&0&0.
}
$$
Hence the bottom row is also exact.
The bimodule $U$ has the property [P2], by Lemmas \ref{transitive}
and \ref{freehasP2}.
Then $U/U\mm$ has property [P1], by Lemma~\ref{PLRtoPL}.
Since $W/W\mm$ is a free left $R$-module then $U/U\mm$ is also
free, by Lemma~\ref{subfree}.
Let $\{b_s\}_{s\in S}$ be a set of elements of $U$ whose residue
classes form a basis of the free left $R$-module $U/U\mm$.
We want to show that this set is a basis of the $R$-$R$-bimodule
$U$.
Since
$$
\bigcap_{i\geq 1}U\mm^i\subseteq\bigcap_{i\geq 1}W\mm^i=\{0\},
$$
then the elements $b_s$, $s\in S$ generate the bimodule $U$,
by Nakayama's lemma.
Assume that
$$
\sum_{s\in S} (r_{s,0}b_s+\sum_{i\geq 2}r_{s,i}b_sm^i)=0,
$$
where all but a finite number of elements $r_{s,i}\in R$ are zero.
It follows from the definition of $b_s$, $s\in S$ that
$r_{s,0}=0$ for any $s\in S$.
Repeating arguments as in the proof of Lemma~\ref{subfree} and
using the fact that $U$ has property [P1'], by Lemma~\ref{bigM},
we get that $r_{s,i}=0$ for any $s\in S$ and $i\geq 2$.
Hence the bimodule $U$ is free.
\end{proof}

\section{The proof of Theorem~\ref{main2}}
\label{proof2}

Suppose that the exact sequence in $\mod A$
\begin{equation} \label{fs}
0\to Z\xrightarrow{\psmatrix{\tilde{f}\\ \tilde{g}}}Z\oplus Y
\xrightarrow{(\tilde{f},-\tilde{h})}Z\to 0
\end{equation}
with $Z$ indecomposable does not split and that $[Z,Z]-[Y,Y]=1$.
Then $\tilde{f}$ is nilpotent and $Y$ is not isomorphic to $Z$.
Furthermore,
$$
[Y,Y]=[Y,Z]=[Z,Y]=[Z,Z]-1,
$$
by Lemmas \ref{seqineq} and \ref{equivsplit} applied to the
sequence \eqref{fs}.
This leads to the following exact sequences induced by \eqref{fs}:
\begin{align}
&0\to\Hom_A(Y,Z)\to\Hom_A(Y,Z\oplus Y)\to\Hom(Y,Z)\to 0,
 \label{4seq1}\\
&0\to\Hom_A(Z,Z)\to\Hom_A(Z,Z\oplus Y)\to\rad(Z,Z)\to 0,
 \label{4seq2}\\
&0\to\Hom_A(Z,Y)\to\Hom_A(Z\oplus Y,Y)\to\Hom(Z,Y)\to 0,
 \label{4seq3}\\
&0\to\Hom_A(Z,Z)\to\Hom_A(Z\oplus Y,Z)\to\rad(Z,Z)\to 0.
 \label{4seq4}
\end{align}

Let $Q$ be the quiver
$$
\xymatrixcolsep{3pc}
\xymatrix@1{
y\ar@/^/[r]^h&z\ar@/^/[l]^g\ar@(ur,dr)[]^f
}
$$
and $\Lambda=kQ/(f^2-hg)$ be the quotient of the path algebra
of $Q$ by the two-sided ideal generated by $f^2-hg$.
We denote by $\varepsilon_y$ and $\varepsilon_z$ the idempotents
corresponding to the vertices $y$ and $z$, respectively.
In particular, $1_\Lambda=\varepsilon_y+\varepsilon_z$.
It is easy to see that
$$
\CB=\{\varepsilon_y,\varepsilon_z,f^{i+1},gf^i,f^ih,gf^ih;\;
i\geq 0\}
$$
is a multiplicative basis of $\Lambda$, that is, $\CB$ is a basis
of the underlying vector space of $\Lambda$ such that $b_1b_2$
belongs to $\CB$ or equals zero, for any $b_1$ and $b_2$ in $\CB$.

Since $(\tilde{f})^2=\tilde{h}\tilde{g}$, we have a canonical
algebra homomorphism
\begin{gather*}
\Phi:\Lambda\to\End_A(Y\oplus Z),\\
\varepsilon_y\mapsto\psmatrix{1_Y&0\\ 0&0},\;
\varepsilon_z\mapsto\psmatrix{0&0\\ 0&1_Z},\;
f\mapsto\psmatrix{0&0\\ 0&\tilde{f}},\;
g\mapsto\psmatrix{0&\tilde{g}\\ 0&0}\;\text{and}\;
h\mapsto\psmatrix{0&0\\ \tilde{h}&0}.
\end{gather*}
This allows us to consider the algebra $E=\End_A(Y\oplus Z)$ as
a $\Lambda$-$\Lambda$-bimodule via
$$
\lambda_1\cdot e\cdot\lambda_2=\Phi(\lambda_1)e\Phi(\lambda_2),
$$
for any $\lambda_1,\lambda_2\in\Lambda$ and
$e\in\End_A(Y\oplus Z)$.
In particular, $\End_A(Y)=\varepsilon_yE\varepsilon_y$.

Let $E'=\psmatrix{\End_A(Y)&\Hom_A(Z,Y)\\ \Hom_A(Y,Z)&\rad(Z,Z)}$.
Then $E'$ is a subbimodule in $E$.
We derive from a direct sum of \eqref{4seq1} and \eqref{4seq2}
the exact sequence
\begin{equation} \label{E1seq}
0\to\varepsilon_zE\xrightarrow{\psmatrix{f\\ g}\cdot}
\varepsilon_zE\oplus\varepsilon_yE\xrightarrow{(f,-h)\cdot}
\varepsilon_zE'\to 0.
\end{equation}

We denote by $R$ the algebra $\varepsilon_y\Lambda\varepsilon_y$
with $1_R=\varepsilon_y$.
The set $\{\varepsilon_y,gf^ih;\;i\geq 0\}$ is a multiplicative
basis of $R$.
Since $(gf^ih)\cdot(gf^jh)=gf^{i+j+2}h$, the algebra $R$ is
commutative and it will be convenient to identify $R$ as
the subalgebra $k[m^2,m^3]$ of the polynomial algebra $k[m]$,
where $gf^ih=m^{i+2}$ for any $i\geq 0$.
Then $\End_A(Y)$ is an $R$-$R$-bimodule.
Let $\{b_1,\ldots,b_s\}$ be a set of generators of the bimodule
$\End_A(Y)$ (for instance we may take a basis of the finite
dimensional vector space $\End_A(Y)$).
Let $\Omega=\Lambda\oplus\bigoplus_{i=1}^s(\Lambda\varepsilon_y
\otimes\varepsilon_y\Lambda)$.
Since $\Lambda\varepsilon_y\otimes\varepsilon_y\Lambda$ is
a projective $\Lambda$-$\Lambda$-bimodule, we may define
the bimodule homomorphism $\Psi:\Omega\to E$,
$$
\Psi(\lambda,\lambda_1\varepsilon_y\otimes\varepsilon_y\lambda'_1,
\ldots,\lambda_s\varepsilon_y\otimes\varepsilon_y\lambda'_s)=
\lambda\cdot 1_E+\sum_{i=1}^s \lambda_i\cdot\psmatrix{b_i&0\\ 0&0}
\cdot\lambda'_i.
$$

\begin{lem}
The homomorphism $\Psi$ is surjective.
\end{lem}

\begin{proof}
Since $b_1,\ldots,b_s$ are generators of $\End_A(Y)$, the latter
is contained in $\im(\Psi)$.
We derive from \eqref{4seq1} the exact squares
$$
\xymatrixcolsep{5pc}
\xymatrix{
\Hom_A(Y,Z)\ar[r]^-{\Hom_A(Y,\tilde{f})}
 \ar[d]_{\Hom_A(Y,\tilde{g})}&\Hom_A(Y,Z)
 \ar[r]^-{\Hom_A(Y,\tilde{g})}\ar[d]_{\Hom_A(Y,\tilde{f})}
 &\End_A(Y)\ar[d]^{\Hom_A(Y,\tilde{h})}\\
\End_A(Y)\ar[r]^-{\Hom_A(Y,\tilde{h})}&\Hom_A(Y,Z)
 \ar[r]^-{\Hom_A(Y,\tilde{f})}&\Hom_A(Y,Z).
}
$$
Consequently, $\Hom_A(Y,Z)=fh\cdot\End_A(Y)+h\cdot\End_A(Y)$
is contained in the $\Lambda$-$\Lambda$-bimodule $\im(\Psi)$.
Similarly, $\im(\Psi)$ contains
$$
\Hom_A(Z,Y)=\End_A(Y)\cdot g+\End_A(Y)\cdot gf.
$$
Let $e\in\End_A(Z)$.
It follows from \eqref{4seq2}
that
$$
e=\mu_1\cdot 1_Z+\tilde{f}e'+\tilde{h}d\qquad\text{and}\qquad
e'=\mu_2\cdot 1_Z+\tilde{f}e''+\tilde{h}d'
$$
for some scalars $\mu_1,\mu_2\in k$ and $A$-homomorphisms $d$,
$d'$, $e'$ and $e''$.
Hence
$$
e=\mu_1\cdot 1_Z+\mu_2\cdot\tilde{f}+\tilde{h}\tilde{g}e''
+\tilde{f}\tilde{h}d'+\tilde{h}d
$$
belongs to
$$
\Psi(\mu_1\cdot\varepsilon_z+\mu_2\cdot f,0,\ldots,0)+h\cdot
\Hom_A(Z,Y)+fh\cdot\Hom_A(Z,Y).
$$
Therefore $\im(\Psi)$ contains $\End_A(Z)$ as well.
\end{proof}

Let $\Lambda'$ denote the subspace of $\Lambda$ generated by
$\CB\setminus\{\varepsilon_z\}$.
Furthermore, let $\Omega'=\Lambda'\oplus\bigoplus_{i=1}^s
(\Lambda\varepsilon_y\otimes\varepsilon_y\Lambda)$.
It is easy to see that $\Lambda'$ is a two-sided ideal
of $\Lambda$ and $\Omega'$ is a $\Lambda$-$\Lambda$-subbimodule
of $\Omega$.

\begin{lem}
The following sequence is exact:
\begin{equation} \label{Omega1seq}
0\to\varepsilon_z\Omega\xrightarrow{\psmatrix{f\\ g}\cdot}
\varepsilon_z\Omega\oplus\varepsilon_y\Omega
\xrightarrow{(f,-h)\cdot}\varepsilon_z\Omega'\to 0.
\end{equation}
\end{lem}

\begin{proof}
$\CB$ induces canonically a basis $\CC$ of the bimodule $\Omega$
such that the set $\CC\cup\{0\}$ is invariant under
left and right multiplications by $f$, $g$ and $h$.
Furthermore, suitable subsets of $\CC$ give bases of the spaces
$\varepsilon_y\Omega$, $\varepsilon_z\Omega$,
$\varepsilon_z\Omega'$, $\Omega\varepsilon_y$,
$\Omega\varepsilon_z$ and $\Omega'\varepsilon_z$.
Now straightforward calculations on these bases are left to
the reader.
\end{proof}

Let $J$ denote the kernel of $\Psi:\Omega\to\End_A(Y\oplus Z)$.

\begin{lem} \label{onesideJ}
The following sequences are exact:
\begin{gather}
\label{J1seq}
0\to\varepsilon_zJ\xrightarrow{\psmatrix{g\\ gf}\cdot}
\varepsilon_yJ\oplus\varepsilon_yJ
\xrightarrow{(fh,-h)\cdot}\varepsilon_zJ\to 0,\\
\label{J2seq}
0\to J\varepsilon_z\xrightarrow{\cdot\psmatrix{h\\ fh}}
J\varepsilon_y\oplus J\varepsilon_y
\xrightarrow{\cdot(gf,-g)}J\varepsilon_z\to 0.
\end{gather}
\end{lem}

\begin{proof}
Observe that $\Psi(\Omega')\subseteq E'$ and let
$e$ denote the element $(\varepsilon_z,0,\ldots,0)$ in $\Omega$.
It follows from the commutative diagram with exact rows
$$
\xymatrix{
0\ar[r]&\Omega'\ar[r]\ar[d]_{\Psi'}&\Omega\ar[r]\ar[d]^\Psi
 &k\cdot\ov{e}\ar[r]\ar@{=}[d]&0\\
0\ar[r]&E'\ar[r]&E\ar[r]&k\cdot\ov{1_Z}\ar[r]&0
}
$$
that $J$ is also the kernel of the restriction
$\Psi':\Omega'\to E'$ of $\Psi$.
Applying \eqref{E1seq} and \eqref{Omega1seq} we get the following
commutative diagram with exact columns and exact two bottom rows:
$$
\xymatrix{
&0\ar[d]&&0\ar[d]&&0\ar[d]\\
0\ar[r]&\varepsilon_zJ\ar[rr]^-{\psmatrix{f\\ g}\cdot}\ar[d]
 &&\varepsilon_zJ\oplus\varepsilon_yJ\ar[rr]^-{(f,-h)\cdot}\ar[d]
 &&\varepsilon_zJ\ar[r]\ar[d]&0\\
0\ar[r]&\varepsilon_z\Omega\ar[rr]^-{\psmatrix{f\\ g}\cdot}\ar[d]
 &&\varepsilon_z\Omega\oplus\varepsilon_y\Omega
 \ar[rr]^-{(f,-h)\cdot}\ar[d]&&\varepsilon_z\Omega'\ar[r]\ar[d]&0\\
0\ar[r]&\varepsilon_zE\ar[rr]^-{\psmatrix{f\\ g}\cdot}\ar[d]
 &&\varepsilon_zE\oplus\varepsilon_yE\ar[rr]^-{(f,-h)\cdot}\ar[d]
 &&\varepsilon_zE'\ar[r]\ar[d]&0\\
&0&&0&&0.
}
$$
Hence the upper row is also exact.
Now the exactness of \eqref{J1seq} follows from joining the
exact squares
$$
\xymatrix{
\varepsilon_zJ\ar[r]^{f\cdot}\ar[d]_{g\cdot}&\varepsilon_zJ
 \ar[r]^{g\cdot}\ar[d]^{f\cdot}&\varepsilon_yJ\ar[d]^{h\cdot}\\
\varepsilon_yJ\ar[r]^{h\cdot}&\varepsilon_zJ\ar[r]^{f\cdot}
 &\varepsilon_zJ.
}
$$
By duality, the sequence \eqref{J2seq} is also exact.
\end{proof}

Observe that $I=\varepsilon_yJ\varepsilon_y$ is an
$R$-$R$-subbimodule of $\varepsilon_y\Omega\varepsilon_y=R\oplus
\bigoplus_{i=1}^s(R\otimes R)$.
Hence we get the exact sequence of $R$-$R$-bimodules
\begin{equation} \label{seqIEnd}
0\to I\to R\oplus\bigoplus_{i=1}^s (R\otimes R)\to\End_A(Y)\to 0.
\end{equation}
In particular, $I$ is torsion free as a right $R$-module.
Obviously the bimodule $R\oplus\bigoplus_{i=1}^s (R\otimes R)$
is generated by the elements $e_0=(1,0,\ldots,0)$ and
$$
e_i=(0,\ldots,0,1\otimes 1,0,\ldots,0),\qquad i=1,\ldots,s.
$$
Furthermore, the bimodule $I$ is finitely generated, since
it is a subbimodule of a finitely generated $R$-$R$-bimodule
and the ring $R\otimes R$ is noetherian.

\begin{lem}
The $R$-$R$-bimodule $I$ has property [P2].
\end{lem}

\begin{proof}
Applying Lemma~\ref{onesideJ} we get the following commutative
diagram with exact rows and columns:
$$
\xymatrix{
&0\ar[d]&&0\ar[d]&&0\ar[d]\\
0\ar[r]&\begin{pmatrix}\varepsilon_zJ\varepsilon_z\end{pmatrix}
 \ar[rr]^-{\cdot\psmatrix{h&fh}}\ar[d]_{\psmatrix{g\\ gf}\cdot}
 &&{\begin{pmatrix}\varepsilon_zJ\varepsilon_y&\varepsilon_zJ
 \varepsilon_y\end{pmatrix}}\ar[rr]^-{\cdot\psmatrix{gf\\ -g}}
 \ar[d]^{\psmatrix{g\\ gf}\cdot}&&\begin{pmatrix}\varepsilon_zJ
 \varepsilon_z\end{pmatrix}\ar[r]\ar[d]^{\psmatrix{g\\ gf}\cdot}
 &0\\
0\ar[r]&{\begin{pmatrix}\varepsilon_yJ\varepsilon_z\\
 \varepsilon_yJ\varepsilon_z\end{pmatrix}}\ar[rr]^-{\cdot
 \psmatrix{h&fh}}\ar[d]_{\psmatrix{fh&-h}\cdot}&&
 {\begin{pmatrix}I&I\\ I&I\end{pmatrix}}\ar[rr]^-{\cdot
 \psmatrix{gf\\ -g}}\ar[d]^{\psmatrix{fh&-h}\cdot}&&
 {\begin{pmatrix}\varepsilon_yJ\varepsilon_z\\ \varepsilon_yJ
 \varepsilon_z\end{pmatrix}}\ar[r]\ar[d]^{\psmatrix{fh&-h}\cdot}
 &0\\
0\ar[r]&\begin{pmatrix}\varepsilon_zJ\varepsilon_z\end{pmatrix}
 \ar[rr]^-{\cdot\psmatrix{h&fh}}\ar[d]&&{\begin{pmatrix}
 \varepsilon_zJ\varepsilon_y&\varepsilon_zJ\varepsilon_y
 \end{pmatrix}}\ar[rr]^-{\cdot\psmatrix{gf\\ -g}}\ar[d]&&
 \begin{pmatrix}\varepsilon_zJ\varepsilon_z\end{pmatrix}\ar[r]
 \ar[d]&0\\
&0&&0&&0.
}
$$
Now the claim follows from the fact that any commutative diagram
in a module category
$$
\xymatrix{
&0\ar[d]&0\ar[d]&0\ar[d]\\
0\ar[r]&B\ar[r]^{\alpha_1}\ar[d]_{\alpha_2}&C_1\ar[r]^{\beta_1}
 \ar[d]^{\gamma_2}&B\ar[r]\ar[d]^{\alpha_2}&0\\
0\ar[r]&C_2\ar[r]^{\gamma_1}\ar[d]_{\beta_2}&D\ar[r]^{\delta_1}
 \ar[d]^{\delta_2}&C_2\ar[r]\ar[d]^{\beta_2}&0\\
0\ar[r]&B\ar[r]^{\alpha_1}\ar[d]&C_1\ar[r]^{\beta_1}\ar[d]
 &B\ar[r]\ar[d]&0\\
&0&0&0
}
$$
with exact rows and columns induces the exact sequence
$$
D\xrightarrow{\gamma_2\alpha_1\beta_2\delta_1}D\xrightarrow{
\psmatrix{\gamma_1\delta_1\\ \gamma_2\delta_2}}D\oplus D.
$$
\end{proof}

Applying Proposition~\ref{freeres} we get an exact sequence of
finitely generated $R$-$R$-bimodules
$$
0\to U\to W\to I\to 0,
$$
where the bimodules $U$ and $W$ are free.
Let $\{u_1,\ldots,u_p\}$ be a basis of the bimodule $U$ and
$\{w_1,\ldots,w_q\}$ be a basis of $W$.
Since the endomorphism $\tilde{f}$ is nilpotent, then
$\tilde{g}(\tilde{f})^{t-2}\tilde{h}=0$ for some $t\geq 2$.
Hence $m^t$ is an annihilator of the left $R$-module $\End_A(Y)$.
Consequently, $I$ contains $m^t(R\oplus\bigoplus_{i=1}^s
(R\otimes R))$, by \eqref{seqIEnd}.
Let $z_i$ be an element of $W$ such that its image is equal to
$m^te_i$ for $i=1,\ldots,s$.

From now on, we shall consider the $R$-$R$-bimodules as left
modules over the algebra $R'=R\otimes R=k[m^2,m^3,n^2,n^3]$,
where the right multiplications by $m^2$ and $m^3$ are replaced
by the multiplications by $n^2$ and $n^3$, respectively.
The algebra $R'$ is contained in $k[m,n]$ and the field $k(m,n)$
of fractions of $k[m,n]$ is also the field of fraction of $R'$.
The tensor product of the $R'$-module monomorphisms
$$
(R\oplus\bigoplus_{i=1}^s (R\otimes R))\xrightarrow{m^t\cdot}
I\to(R\oplus\bigoplus_{i=1}^s (R\otimes R))
$$
by the flat $R'$-module $k(m,n)$ leads to the monomorphisms
$$
k(m,n)^s\to I\otimes_{R'}k(m,n)\to
k(m,n)^s.
$$
Since the composition is an isomorphism, each of the above maps
is an isomorphism of vector spaces over the field $k(m,n)$.
It follows from the exact sequence
$$
0\to U\otimes_{R'}k(m,n)\to W\otimes_{R'}
k(m,n)\to I\otimes_{R'}k(m,n)\to 0
$$
that $p+s=q$.
Let $D$ be the $q\times q$-matrix with coefficients in $R'$ such
that its rows represent the elements
$u_1,\ldots,u_p,z_1,\ldots,z_s$ in the basis $w_1,\ldots,w_q$,
that is,
$$
\bsmatrix{u_1\\ \vdots\\ z_s}=D\bsmatrix{w_1\\ \vdots\\ w_q}.
$$

\begin{lem} \label{detD}
$\det(D)=m^j$ for some integer $j$.
\end{lem}

\begin{proof}
Let $z_0$ be an element of $W$ whose image in $I$ is equal to
$m^te_0$.
Then $(m^2-n^2)z_0$ belongs to $U$ and hence
$$
(m^2-n^2)z_0=r'_1u_1+\ldots+r'_pu_p
$$
for some elements $r'_i$ in $R'$.
This implies that
$$
z_0=\frac{1}{m^2-n^2}\cdot v_2\cdot\bsmatrix{u_1\\ \vdots\\ z_s},
\qquad\text{where }v_2=[r'_1,\ldots,r'_p,0\ldots,0].
$$
Since the image of $m^tw_i$ belongs to the $R'$-submodule
in $I$ generated by $m^te_0,m^te_1,\ldots,m^te_s$, then
$m^tw_i$ belongs to the $R'$-submodule in
$W$ generated by $u_1,\ldots,u_p,z_1,\ldots,z_s$ and $z_0$,
for $i=1,\ldots,q$.
Consequently,
there are a $q\times q$-matrix $B$ with coefficients in $R'$
and elements $r_1,\ldots,r_q$ in $R'$ such that
$$
\bsmatrix{m^tw_1\\ \vdots\\ m^tw_q}=B\bsmatrix{u_1\\ \vdots\\ z_s}
+\bsmatrix{r_1\\ \vdots\\ r_q}\cdot z_0.
$$
Observe that
$$
\bsmatrix{m^tw_1\\ \vdots\\ m^tw_q}=C\bsmatrix{u_1\\ \vdots\\ z_s},
$$
where
$$
C=B+\frac{1}{m^2-n^2}\cdot\bsmatrix{r_1\\ \vdots\\ r_q}\cdot v_2
$$
is a $q\times q$-matrix with coefficients in the field $k(m,n)$.
Consequently, $C\cdot D=m^t\cdot I_q$, where $I_q$ denotes the
identity matrix.
Hence
$$
\det(C)\cdot\det(D)=m^{tq}.
$$
Let $D_i$ be the matrix obtained from $B$ by replacing the $i$-th
row by $v_2$, for $i=1,\ldots,q$.
Applying elementary properties of determinants we get
$$
\det(C)=\det(B)+\frac{r_1}{m^2-n^2}\det(D_1)+\ldots
+\frac{r_q}{m^2-n^2}\det(D_q).
$$
Therefore $(m^2-n^2)\cdot\det(C)$ is an element of
$k[m^2,m^3,n^2,n^3]$ and
$\det(D)$ is a divisor of $(m^2-n^2)\cdot m^{tq}$ in the algebra
$k[m^2,m^3,n^2,n^3]$.
Replacing $(m^2-n^2)$ by $(m^3-n^3)$, we get that
$\det(D)$ is also a divisor of $(m^3-n^3)\cdot m^{tq}$.
Since $k[m^2,m^3,n^2,n^3]$ is contained in the unique
factorization domain $k[m,n]$ and the polynomials
$$
\frac{m^2-n^2}{m-n}=m+n \qquad\text{and}\qquad
\frac{m^3-n^3}{m-n}=m^2+mn+n^2
$$
are coprime, then $\det(D)$ is a divisor of $(m-n)m^{tq}$.
Therefore $\det(D)=m^j$ for some integer $j$, as $(m-n)m^j$ does
not belong to $k[m^2,m^3,n^2,n^3]$.
\end{proof}

Applying Lemma~\ref{detD} we get that the coefficients of
the matrix $m^j\cdot D^{-1}$ belong to $k[m^2,m^3,n^2,n^3]$.
Hence $m^jx$ belongs to the $R'$-submodule of $W$ generated by
$u_1,\ldots,u_p,z_1,\ldots,z_s$, for any $x\in W$.
Therefore $m^jy$ belongs to the $R'$-submodule of $I$
generated by $m^te_1,\ldots,m^te_s$, for any $y\in I$.
Taking $y=m^t\cdot e_0$ we obtain a contradiction.
This finishes the proof of Theorem~\ref{main2}.

\section{Corollaries and remarks}
\label{remarks}

\textbf{\ref{remarks}.1.}
Theorem~\ref{main1} can be generalized to other varieties.
We give here two examples.

Let $Q=(Q_0,Q_1,s,e)$ be a finite quiver, where $Q_0$ is the set
of vertices, $Q_1$ is the set of arrows and $s,e:Q_1\to Q_0$ are
functions such that any arrow $\alpha$ in $Q_1$ has the starting
vertex $s(\alpha)$ and the ending vertex $e(\alpha)$.
Let $\dd=(d_i)_{i\in Q_0}$ be a sequence of positive integers.
Furthermore, we denote by $\BM_{d'\times d''}(k)$ the space of
$d'\times d''$-matrices with coefficients in $k$, for any positive
integers $d'$ and $d''$.
Then the group $\GL_\dd(k)=\prod_{i\in Q_0}\GL_{d_i}(k)$ acts on
the affine space $\rep_Q^\dd(k)=\prod_{\alpha\in Q_1}
\BM_{d_{e(\alpha)}\times d_{s(\alpha)}}(k)$ by conjugations
$$
(g_i)_{i\in Q_0}\star(m_\alpha)_{\alpha\in Q_1}=(g_{e(\alpha)}
m_\alpha g_{s(\alpha)}^{-1})_{\alpha\in Q_1}.
$$
Let $d=\sum_{i\in Q_0}d_i$ and $kQ$ denote the path algebra of $Q$.
Then there is a fibre bundle
$$
\CC_\dd\to\left(\GL_d(k)/\GL_\dd(k)\right)
$$
with typical fiber $\rep_Q^\dd(k)$, where $\CC_\dd$ is a connected
component of $\mod_{kQ}^d(k)$ (see \cite{Bgeo}).
Consequently, Theorem~\ref{main1} remains true if we take the
$\GL_\dd(k)$-variety $\rep_Q^\dd(k)$ instead of $\mod_A^d(k)$.

Let $P_1,\ldots,P_t$ be parabolic subgroups of $G=\GL_d(k)$.
We consider the projective variety
$$
X=G/P_1\times\ldots\times G/P_t
$$
equipped with the diagonal action of $G$.
Applying arguments used in \cite[\S 2]{BZ2} we get
a $G$-equivariant principal $H$-bundle $\CU\to X$, where
$\GL_\dd(k)=G\times H$ and $\CU$ is a $\GL_\dd(k)$-invariant open
subset of $\rep_Q^\dd(k)$, for some quiver $Q$ and sequence $\dd$.
Thus Theorem~\ref{main1} is still true if we replace the module
variety by the $G$-variety $X$.
\bigskip

\noindent\textbf{\ref{remarks}.2.}
Let $Q$ be the Kronecker quiver.
Then the orbit closures in $\rep_Q^\dd(k)$ are regular in
codimension one, even if they contain infinitely many orbits
(\cite{BB}).
It is an interesting question whether the orbit closures are
regular in codimension one for the other extended Dynkin quivers.
\bigskip

\noindent\textbf{\ref{remarks}.3.}
We say that two exact sequences in $\mod A$
$$
\sigma_l: 0\to Z_l\xrightarrow{f_l}Z_l\oplus M\xrightarrow{g_l}N
\to 0,\qquad l=1,2,
$$
with $f_l$ in $\rad(Z_l,Z_l\oplus M)$ are equivalent if there is
a commutative diagram in $\mod A$
$$
\xymatrix{
0\ar[r]&Z_1\ar[r]^-{f_1}\ar[d]_i&Z_1\oplus M\ar[r]^-{g_1}\ar[d]^j&
 N\ar[r]\ar@{=}[d]&0\\
0\ar[r]&Z_2\ar[r]^-{f_2}&Z_2\oplus M\ar[r]^-{g_2}&N\ar[r]&0
}
$$
for some isomorphisms $i$, $j$.
In particular, $Z_1$ is isomorphic to $Z_2$.

\begin{cor}
Let $M$ and $N$ be points in $\mod_A^d(k)$ such that $N$ belongs
to $\ov{\CO}_M$ and $\dim\CO_M-\dim\CO_N=1$.
Then there is a unique, up to an equivalence, exact sequence
in $\mod A$
\begin{equation} \label{onceagain}
0\to Z\xrightarrow{f}Z\oplus M\xrightarrow{g} N\to 0
\end{equation}
with a radical morphism $f$.
Furthermore, the module $Z$ is indecomposable.
\end{cor}

\begin{proof}
Applying Lemma~\ref{Zind} we get the exact sequence
\eqref{onceagain} with $f$ radical and $Z$ indecomposable.
Let
$$
\sigma': 0\to Z'\xrightarrow{f'}Z'\oplus M\xrightarrow{g'}N\to 0
$$
be an exact sequence in $\mod A$ with $f'$ in
$\rad(Z',Z'\oplus M)$.
By Lemma~\ref{dimhomend} and Proposition~\ref{longprop},
$[Z\oplus M,M]=[Z\oplus M,N]$.
Applying Lemma~\ref{seqineq} for $\sigma'$ and $X=Z\oplus M$ we get
that $g$ factors through $g'$.
This leads to a commutative diagram in $\mod A$
$$
\xymatrix{
0\ar[r]&Z\ar[r]^-f\ar[d]_i&Z\oplus M\ar[r]^g\ar[d]^j&
 N\ar[r]\ar@{=}[d]&0\\
0\ar[r]&Z'\ar[r]^-{f'}&Z'\oplus M\ar[r]^-{g'}&N\ar[r]&0
}
$$
for some homomorphisms $i$ and $j$.
We conclude from Lemma~\ref{equivsplit} that the induced exact
sequence
$$
0\to Z\xrightarrow{\psmatrix{f\\ i}}(Z\oplus M)\oplus Z'
\xrightarrow{(j,-f')}(Z'\oplus M)\to 0
$$
splits.
Hence $i$ is a section and $j$ is a retraction, as $f$ and $f'$
are radical homomorphisms.
This implies that $Z$ is a direct summand of $Z'$ as well as
$Z'\oplus M$ is a direct summand of $Z\oplus M$.
Consequently, $Z$ is isomorphic to $Z'$, and $i$ and $j$ are
isomorphisms.
\end{proof}
\bigskip

\noindent\textbf{\ref{remarks}.4.}
The bound for the difference of dimensions of $\End_A(Z)$ and
$\End_A(Y)$ given in Theorem~\ref{main2} is sharp.
We recall here the example given in \cite[5.4]{ARZ}.
Let $A=k[\alpha,\beta]/(\alpha^2,\beta^2)$ and $Y$ and $Z$ be
modules in $\mod_A^4(k)$ such that
$$
Y(\alpha)=Z(\alpha)=\psmatrix{
0&0&0&0\\
0&0&0&0\\
1&0&0&0\\
0&1&0&0
},\quad Y(\beta)=\psmatrix{
0&0&0&0\\
1&0&0&0\\
0&0&0&0\\
0&0&1&0
},\quad Z(\beta)=\psmatrix{
0&0&0&0\\
0&0&0&0\\
0&0&0&0\\
1&0&0&0
},
$$
and set
$$
\tilde{f}=\psmatrix{
0&0&0&0\\
1&0&0&0\\
0&1&0&0\\
0&0&1&0
},\quad \tilde{g}=\psmatrix{
0&0&0&0\\
0&1&0&0\\
1&0&0&0\\
0&0&0&1
},\quad \tilde{h}=\psmatrix{
1&0&0&0\\
0&0&0&0\\
0&0&1&0\\
0&1&0&0
}.
$$
Then the module $Z$ is indecomposable, the sequence
$$
0\to Z\xrightarrow{\psmatrix{\tilde{f}\\ \tilde{g}}}Z\oplus Y
\xrightarrow{(\tilde{f},-\tilde{h})}Z\to 0
$$
is exact and $\dim_k\End_A(Z)-\dim_k\End_A(Y)=2$.
\bigskip

\noindent\textbf{\ref{remarks}.5.}
The properties [P1] and [P2] are strongly related to free
resolutions of modules over $R=k[m^2,m^3]$ and
$R'=k[m^2,m^3,n^2,n^3]$, respectively.

\begin{cor} \label{finresbimod}
Let $M$ be an $R'$-module.
Then the following conditions are equivalent:
\begin{enumerate}
\item[(1)] $M$ has property [P2];
\item[(2)] there is a free resolution $0\to F_2\to F_1\to F_0\to
M\to 0$ of $M$;
\item[(3)] $M$ has a free resolution of finite length.
\end{enumerate}
\end{cor}

\begin{proof}
Obviously \textit{(2)} implies \textit{(3)}.
Furthermore, \textit{(3)} implies \textit{(1)}, by
Lemmas \ref{transitive} and \ref{freehasP2}.
Assume now that $M$ has property [P2].
We take an exact sequence of $R'$-modules
$$
0\to M'\to F_0\to M\to 0,
$$
where the module $F_0$ is free.
The module $M'$ has property [P2], by Lemmas \ref{transitive} and
\ref{freehasP2}.
Since $M'$ is a torsion free bimodule, then \textit{(2)} follows
from Proposition~\ref{freeres}.
\end{proof}

\begin{cor} \label{finresmod}
Let $M$ be an $R$-module.
Then the following conditions are equivalent:
\begin{enumerate}
\item[(1)] $M$ has property [P1];
\item[(2)] there is a free resolution $0\to F_1\to F_0\to M\to 0$
of $M$;
\item[(3)] $M$ has a free resolution of finite length.
\end{enumerate}
\end{cor}

\begin{proof}
The proof is similar to the previous one.
We have to replace Proposition~\ref{freeres} by
Lemma~\ref{subfree}.
Furthermore, one can repeat appropriate arguments to get
versions of Lemmas \ref{transitive} and \ref{freehasP2}
for $R$-modules.
\end{proof}

Observe that $N=\begin{pmatrix}M\\ M\end{pmatrix}$ is
a $k[\xi]/(\xi^2)$-module for any $R$-module $M$ and
$N'=\begin{pmatrix}M'&M'\\ M'&M'\end{pmatrix}$ is
a $k[\xi,\eta]/(\xi^2,\eta^2)$-module for any $R'$-module $M'$,
where the residue classes of $\xi$ and $\eta$ denote the
multiplications $\psmatrix{m^3&-m^2\\ m^4&-m^3}\cdot$
and $\cdot\psmatrix{n^3&n^4\\ -n^2&-n^3}$, respectively.
Since the algebras $k[\xi]/(\xi^2)$ and
$k[\xi,\eta]/(\xi^2,\eta^2)$ are local and Frobenius, the
free modules over them coincide with the projective modules
and with the injective ones.
One can prove that $M$ has property [P1] if and only if
$N$ is a free $k[\xi]/(\xi^2)$-module, and $M'$ has property [P2]
if and only if $N'$ is a free $k[\xi,\eta]/(\xi^2,\eta^2)$-module.

\end{document}